\documentclass[a4paper,english,pdftex]{bourbaki}

\usepackage[colorlinks=true,linkcolor=blue,citecolor=blue]{hyperref}
\usepackage{amsmath}
\usepackage{amssymb}
\usepackage{mathrsfs}
\usepackage[all]{xy}
\usepackage{graphicx}
\usepackage{latexsym}
\usepackage{amsthm}
\usepackage{amscd}
\usepackage{mathrsfs}

\def\qbar{{\bar{\bold Q}}}
\def\qpbar{{\bar{\bold Q}_p}}

\def\fpbar{{\bar{\bold F}_p}}

\DeclareMathOperator\aut{Aut}
\DeclareMathOperator\bel{Bel}
\DeclareMathOperator\gal{Gal}

\title{On the Belyi degree(s) of a curve defined over a number field}
\author{Leonardo Zapponi}
\date{\today}

\begin{document}

\maketitle

\tableofcontents

\section*{Introduction} The celebrated theorem of Belyi in~\cite{Belyi} asserts that any smooth projective curve defined over $\qbar$ can be realized as a Belyi cover (or Belyi map), that is a cover $f:X\to\bold P^1$ unramified outside three points, the converse being true by a rigidity criterion of Weil, cf.~\cite{Weil}. This result can be considered as an arithmetic uniformization of curves. It is then natural to investigate on the properties of the Belyi covers $f$ defined over a fixed curve $X$ and, more particularly, on how the geometry of $f$ is related to the arithmetic of $X$. From a geometric point of view, the behaviour of $f$ only depends on the $\qbar$-isomorphism class of $X$. In particular, if no restrictions are made on the field of definition of $f$ (for example) then we are forced to consider those arithmetic invariants which are stable under $\qbar$-isomorphism.

A first natural question is to determine the minimal degree of a Belyi cover on $X$. We refer to this integer as the (absolute) Belyi degree of $X$. This invariant was first introduced and studied by R. Li{\c{t}}canu in~\cite{Litcanu} (see also~\cite{Litcanu2}) who proved that the Belyi degree somehow behaves like a height function. More precisely, Grothendieck's theory of Dessins d'enfants implies that there exist finitely many isomorphism classes of Belyi covers of bounded degree (see~\cite{Dessin} for an introduction to this subject). As a direct consequence, there exist finitely many $\qbar$-isomorphism classes of curves with bounded Belyi degree, cf. Proposition~\ref{finite1}. It is then interesting to find some bounds for this invariant. As it is done in~\cite{Litcanu}, the combinatorial approach leads to an upper bound obtained by counting particular classes of graphs. On the other hand, a slight modification of the techniques used in~{\it loc. cit.} for the study of the Belyi degree of an algebraic number leads to an upper bound with a more arithmetical flavour, depending on the (usual) height of the branch locus of a particular cover $X\to\bold P^1$ (of degree bounded by the genus, for example). The Rieman-Hurwitz formula gives a lower bound for the Belyi degree, only depending on the genus of the curve, but this result is not really interesting. In Theorem~\ref{sharp}, we give a new lower bound related to the reduction behaviour of the curve. More precisely, we define the notion of stable (rational) prime of bad reduction, which only depends on the $\qbar$-isomorphism class of the curve. We then prove that the Belyi degree is greater than or equal to the greatest stable prime of bad reduction (if it exists). This result is a direct consequence of Beckmann's results in~\cite{Beckmann}. At a first glance, it would seem that this bound is quite rough, but we then show that for any positive integer $g$ and any prime number $p>2g+1$, there exists a genus $g$ curve of Belyi degree $p$ having $p$ as stable prime of bad reduction. In particular, if no additional assumptions are made, the bound obtained here can be considered as optimal. The proof of this result is the only place in the paper where dessins d'enfants explicitely appear. No particular knowledge of this theory is needed, we just use them in order to prove that a certain cover exists (which, in itself, is a highly non-trivial fact).

As it was already pointed out by R. Li{\c{t}}canu, with the present definition of the Belyi degree we canot expect a finiteness result as Proposition~\ref{finite1} when restricting to $K$-isomorphism classes (rather than $\qbar$-isomorphism classes), where $K$ is a number field. This is the main difference with (what is expected to hold for) a height function. In order to overcome this problem, a first idea is to define the relative Belyi degree of a curve $X$ defined over a (fixed) number field $K$ as the minimal degree of a Belyi map $X\to\bold P^1$ defined over $K$. Unfortunately, in despite of this restriction, there still can exist infinitely many $K$-isomorphism classes of curves of bounded (relative) Belyi degree. This is related to the fact that a Belyi map can have non-trivial $K$-twists. Since this situation only occurs if the automorphism group of the cover is non-trivial, we finally define the relative Belyi degree by only considering the Belyi maps defined over $K$ and with trivial automorphism group. With this definition, the finiteness result applies for $K$-isomorphism classes. Moreover, the relative Belyi degree is greater than or equal to the absolute Belyi degree and it behaves well under finite extension of the number field $K$ (it is non-increasing and stationnary as soon as the extension contains a certain number field $L$). There just remains to prove that any curve $X$ defined over $K$ admits such a Belyi map. This is done in Proposition~\ref{existence}. The proof is quite long, but essentially elementary and based on linear algebra arguments. After writing this paper, we found an earlier and essentially identical proof of this result in~\cite{Couveignes}, we nevertheless include it in this preprint version and will omit it if this paper will be published. If the genus of the curve is greater than $1$ then the proof can be simplified but we decided to give an unified and effective proof including the genus $0$ case, since we can now also define the relative Belyi degree of a conic defined over a number field. In fact, the genus zero case is interesting and defenitely not trivial when working with $K$-isomorphism classes. Some illuminating examples can be found in~\cite{Couveignes2,Couveignes3}, which can be considered as being essentially the only papers treating this question. In Theorem\ref{bound1}, we obtain a lower bound for the relative Belyi degree which can be considered as the relative version of Theorem~\ref{bound1}: we show that it is greater than or equal to the greatest prime of bad reduction of the minimal regular model of the curve (as soon as its genus is positive, the case of conics may be treated differently, but Th\'eor\`eme 6 in~\cite{Couveignes3} can be considered as the right analogue). This result once again follows from a slight modification of Beckmann's result, which asserts that if $X\to\bold P^1$ is a Belyi map defined over $K$ with trivial automorphism group then the (minimal model of the) curve $X$ has good reduction at any prime of $K$ not dividing the order of the monodromy group of the cover. This result is well-known to experts but the author could not find a reference for it. Its short proof has therefore been included. In this paper, we restricted to projective curves but everything carries out naturally when working in the affine case (some minor modifications are needed when defining the notion bad reduction). For example, if the curve is (isomorphic to) the projective line minus $4$ points, we recover the (absolute) Belyi degree of an algebraic number, as defined in~\cite{Litcanu}, and its (new) relative version.

\section*{Aknowledgements} This paper was written after an informal aftenoon discussion with W. Goldring around Belyi's Theorem. I would like to thank him for giving me the impulse to write down these results, which were known to me since quite a long time. I am also indebted to S. Wewers for some comments concerning the proof of Theorem~\ref{bound2}.

\section{The absolute Belyi degree of a curve}

\subsection{Definition and first properties} Let $X$ be a smooth projective curve of genus $g$ defined over $\qbar$. Following Belyi's Theorem~\cite{Belyi}, there exists a finite cover $f:X\to\bold P^1$ unramified outside the set $\{\infty,0,1\}$, generally called Belyi map. The {\bf (absolute) Belyi degree} of $X$, which will be denoted by $\deg_B(X,\qbar)$, is the minimal degree of such a cover. If $\bel(X,\qbar)$ denotes the set of Belyi maps $X\to\bold P^1$, we clearly have the identity
$$\deg_B(X,\qbar)=\min\{\deg(f)\,\,|\,\,f\in\bel(X,\qbar)\}.$$
This definition clearly only depends on the $\qbar$-isomorphism class of $X$ and the Belyi degree can thus be considered as a map
$$\mathcal M_g(\qbar)\to\bold N,$$
where $\mathcal M_g$ is the (coarse) moduli space of genus $g$ curves. The following result directly follows from the definitions, see also~\cite{Litcanu,Litcanu2} for further details.

\begin{prop}\label{finite1}There exist finitely many $\qbar$-isomorphism classes of curves of bounded Belyi degree.
\end{prop}

\begin{exem} The projective line is the unique curve of Belyi degree less than or equal to $2$. There is a unique isomorphism class of curves of Belyi degree $3$, corresponding to the ellitpic curves with $j$-invariant $0$. There are two isomorphism classes of elliptic curves with Belyi degree $4$, corresponding to $j=1728$ and $j=\frac{207646}{6561}$.\end{exem}

\subsection{Stable primes of bad reduction, semi-stable models} The curve $X$ being as above, let $p$ be a prime number and fix an injection $\qbar\hookrightarrow\qpbar$. By base change, we can then view $X$ as a curve over $\qpbar$ and the Semi-Stable Reduction Theorem~\cite{Deligne-Mumford} asserts that there exists a stable model $\mathcal X$ of $X$ defined over the ring of integers $R$ of a $p$-adic field $K\subset\qpbar$. The $R$-curve $\mathcal X$ is unique up to $R$-isomorphism and, as above, it only depends on the $\qpbar$-isomorphism class of $X$. We then say that $X$ has {\bf (potentially) good reduction} if the reduced curve $\bar{\mathcal X}$ is smooth, otherwise $X$ has {\bf bad reduction}. We say that $p$ is a {\bf stable prime of bad reduction} if there exists an injection $\qbar\hookrightarrow\qpbar$ such that (the stable model of) $X$ has bad reduction. The set of stable primes of bad reduction only depend on the $\qbar$-isomorphism class of $X$. In terms of moduli spaces, the existence of the semi-stable model can be viewed as a specialization map
$$\mathcal M_g(\qpbar)\to\bar{\mathcal M}_g(\fpbar),$$
where $\bar{\mathcal M}_g$ denotes the Deligne-Mumford compatification of $\mathcal M_g$, cf. {\it loc. cit}.

\subsection{A lower bound}

\begin{theo}\label{bound1} For any curve $X$ defined over $\qbar$, we have the inequality
$$\deg_B(X,\qbar)\geq p,$$
where $p$ is the greatest stable prime of bad reduction.
\end{theo}

\begin{proof} Let $f:X\to\bold P^1$ be a Belyi map defined over $\qbar$ and denote by $G$ its monodromy group. Let $p$ be a prime not dividing the order of $G$. Following Proposition 5.3 in~\cite{Beckmann}, there exists a model of the cover defined over the ring of integers of a number field $K$ having good reduction at any prime lying above $p$ (the original proof only deals with the Galois case but the result is true in full generality). This implies that for any injection $\qbar\hookrightarrow\qpbar$, the stable model of curve $X$ has good reduction and thus $p$ cannot be a stable prime of bad reduction. In other words, the stable primes of bad reduction of $X$ divide the order of $G$. The result then follows from the fact that the group $G$ can be realized as a subgroup of the symmetric group $S_n$, where $n=\deg(f)$, so that any prime divisor of the order of $G$ is less than or equal to $n$.
\end{proof}

\subsection{Sharpness of the bound}

\begin{theo}\label{sharp} For any positive integer $g$ and for any prime number $p\geq 2g+3$, there exists a curve of genus $g$ defined over $\qbar$ of (absolute) Belyi degree $p$ and having $p$ as stable prime of bad reduction.
\end{theo}

\begin{proof} Given the integers $g$ and $p$, consider a Belyi map $f:\bold P^1\to\bold P^1$ associated to the dessin d'enfant of Figure 1.
\vskip.4cm
\begin{center}
\includegraphics[scale=.6]{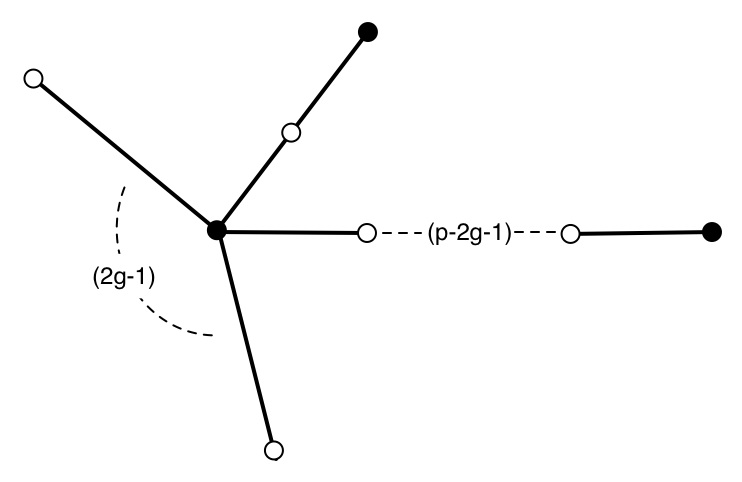}\\
Figure 1.
\end{center}
\vskip.4cm
The degree $p$ cover $f$ has the following ramification datum:
\begin{enumerate}
\item There is a unique point lying above $\infty$.
\item There are $\frac{p+3}2-g$ points above $0$, one of them has ramification index $2g+1$, two are unramified and the remaining have ramification index $2$.
\item There are $\frac{p-1}2+g$ points above $1$, $2g-1$ of them are unramified and the others have ramification index $2$.
\end{enumerate}
There exist exactly $2g+2$ unramified points in the union of the ramified fibers of $f$, three of them above $0$ and the remaining above $1$. Let $X$ be the smooth hyperelliptic genus $g$ curve obtained as double cover of the projective line ramified above these points. Since the cover $f$ may be defined over $\qbar$, the curve $X$ is itself defined over $\qbar$. From an explicit point of view,  the cover $f$ corresponds to an inclusion $\qbar(x)\hookrightarrow\qbar(t)$. In this case, the field $\qbar(t,\sqrt{x(1-x)})$ is the field of rational functions of $X$. By construction, the degree $2p$ cover $\tilde f:X\to\bold P^1$ obtained by composition of $f$ with the canonical projection is a Belyi map with the following ramification datum:
\begin{enumerate}
\item There are two points above $\infty$, with ramification index $p$.
\item There are $p-2g$ points above $0$, one of them has ramification index $4g+2$ while the others have ramification index $2$.
\item There are $p$ points above $1$, all of them having ramification index $2$.
\end{enumerate}
Now, the cover $\tilde f$ factors as
$$X\stackrel h\longrightarrow\bold P^1\stackrel\rho\longrightarrow\bold P^1,$$
where the double cover $\rho$ only ramifies above $0$ and $1$. It then follows from Abhyankar's Lemma that the degree $p$ cover $h$ is unramified outside three points; it then defines a Belyi map with the following ramification datum:
\begin{enumerate}
\item There is a unique point above $\infty$, with ramification index $p$.
\item There is a unique point above $0$, with ramification index $p$.
\item There are $p-2g$ points above $1$, one of them has ramification index $2g+1$ and all the others are unramified.
\end{enumerate}
We finally study the reduction behaviour of $X$. Fix an embedding $\qbar\hookrightarrow\qpbar$, so that $X$ can be considered as a curve over $\qpbar$. The results in~\cite{Zapponi1} give an explicit description of the semi-stable model $\mathcal C\to\mathcal D$ of $f$ (more precisely, the minimal semi-stable model which separates the ramified locus, cf. {\it loc. cit.} for a precise definition). For our purpose, we just need the following facts:
\begin{enumerate}
\item The reduced curve $\bar{\mathcal C}$ (resp. $\bar{\mathcal D}$) is the union of three projective lines $C_\infty,C_0$ and $C_1$ (resp. $D_\infty,D_0$ and $D_1$) such that $C_\infty$ meets $C_0$ and $C_1$ in two distinc points (resp. $D_\infty$ meets $D_0$ and $D_1$ in two distinct points), these being the only singularities, cf. the figure below.
\item The point $\lambda\in\{\infty,0,1\}$ specializes in (the smooth locus of) $D_\lambda$ and the fiber above $\lambda$ specializes in (the smooth locus of) $C_\lambda$.
\end{enumerate}
\vskip.4cm
\begin{center}
\includegraphics[scale=.6]{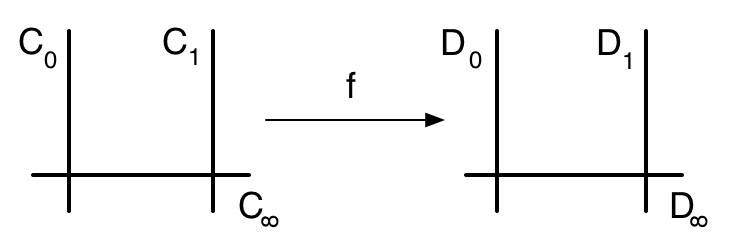}\\
Figure 2.
\end{center}
\vskip.4cm
By working with Weierstrass equations  (see for example~\cite{Liu}), it then easily follows that the stable model of $X$ is the union of an elliptic curve and an hyperelliptic curve of genus $g-1$ meeting at a single point, which implies that $p$ is a stable prime of bad reduction.
\end{proof}

\section{The relative Belyi degree of a curve}

\subsection{Definition and first properties} Let $K$ be a number field. As stated before, given a smooth projective curve $X$ defined over $K$, the existence of a Belyi map $f:X\to\bold P^1$ of given degree only depends on the $\qbar$-isomorphism class of $X$. With the present definition of the (absolute) Belyi degree, we cannot expect a finiteness result similar to Proposition~\ref{finite1} concerning $K$-isomorphism classes of curves. Following the notation of the first section, let $\bel(X,K)$ be the subset of $\bel(X,\qbar)$ consisting of Belyi maps defined over $K$. For any $f\in\bel(X,K)$, denote by
$$\aut_\qbar(f)=\{\sigma\in\aut_\qbar(X)\,\,|\,\, f\circ\sigma=f\}$$
its $\qbar$-automorphism group and consider the set
$$\bel^0(X,K)=\{f\in\bel(X,K)\,\,|\,\,\aut(f)=1\}$$
consisting of Belyi maps with trivial automorphism group. The following classical and elementary result asserts that any element of $\bel^0(X,K)$ separates the $K$-isomorphism class defined by $X$ inside its $\qbar$-isomorphism class.

\begin{lemm}\label{lemm1} Let $f\in\bel^0(X,K)$ and suppose that $\sigma:Y\to X$ a $\qbar$-isomorphism, where $Y$ is defined over $K$. Then the induced cover $g=f\circ\sigma:Y\to\bold P^1$ is defined over $K$ if and only if $\sigma$ is defined over $K$.
\end{lemm}

\begin{proof} One implication is trivial, so assume that the Belyi map $g$ is defined over $K$. For any element $\tau\in\gal(\qbar/K)$, we obtain an element
$$\phi_\tau={}^\tau\sigma\circ\sigma^{-1}\in\aut_\qbar(X).$$
We then have the identities
$$f\circ\phi_\tau={}^\tau f\circ{}^\tau\sigma\circ\sigma^{-1}={}^\tau g\circ\sigma^{-1}=g\circ\sigma^{-1}=f$$
which imply that $\phi_\tau$ actually belongs to $\aut_\qbar(f)$, so it is the identity. We finally obtain the relation ${}^\tau\sigma=\sigma$, so that $\sigma$ is defined over $K$.
\end{proof}

In order to define the relative Belyi degree, we must prove that any curve admits a Belyi cover with trivial automorphism group.

\begin{prop}\label{existence} For any smooth projective curve $X$ defined over a number field $K$, the set $\bel^0(X,K)$ is not empty.
\end{prop}

\begin{proof} The proof is essentially the same as in~\cite{Couveignes}, we nevertheless decided to include it in the paper. Let $S_0\subset X(\qbar)$ be a finite subset stable under the action of $\mbox{Gal}(\qbar/K)$ and such that the group $G=\aut_\qbar(X-S_0)$ is finite. If $g$ denotes the genus of $X$ and $n$ the cardinality of $S$, this last condition is equivalent to the inequality $2g-2+n>0$. Viewing $S_0$ as a divisor on $X$, which is ample, fix an integer $r$ such that $D=rS_0$ is very ample and consider the invertible sheaf $\mathcal L=\mathcal O_X(D)$. Let $L/K$ be a finite extension containing $\mu_{|G|}$ and the field of definition of $G$. Denote by $X'=X\otimes_KL$ the curve obtained from $X$ by base change. The sheaf $\mathcal L$ then defines a very ample invertible sheaf $\mathcal L'$ on $X'$ such that, setting $V=H^0(X,\mathcal L)$ and $V'=H^0(X',\mathcal L')$, we have the identity $V'=V\otimes_KL$. In the following, we consider $V$ as a sub-$K$-vector space of $V'$. Now, the group $G$ naturally acts $L$-linearly on $V'$. The very ampleness of $\mathcal L'$ implies that the action is faithfull. Given an element $\sigma\in G$, its eigenvalues all belong to $\mu_{|G|}$. For any $\sigma\in G-\{1\}$ and any $\zeta\in\mu_{|G|}$, let $V'(\sigma,\zeta)\subset V'$ be the eigenspace of $\sigma$ corresponding to $\zeta$ (which may be trivial) and consider the sub-$K$-vector space $V(\sigma,\zeta)=V'(\sigma,\zeta)\cap V$. Since $\sigma$ is a $L$-automorphism of $L(X)$, it follows that $V(\sigma,\zeta)$ is properly contained in $V$. Indeed, if $\zeta=1$, the identity $V(\sigma,\zeta)=V$ would lead to $V'(\sigma,\zeta)=V'$, which is excluded by the faithfulness of the action of $G$. Similarly, for $\zeta\neq 1$, the identity $V(\sigma,\zeta)=V$ would give $V'(\sigma,\zeta)=V'$, which is impossible since $L\subset V'(\sigma,1)$. In particular, since $K$ is infinite and $G$ is finite, the $K$-vector space $V$ cannot be the union of the $V(\sigma,\zeta)$. Let $f\in V$ be an element not belonging to this union and denote by $S_1-S_\infty$ its divisor, where $S_1$ and $S_\infty$ are positive and their support have trivial intersection. We can moreover assume that $S_\infty=rS_0$ (the argument is the same as above). By construction, the divisor $S_1$ is stable under the action of $\mbox{Gal}(\qbar/K)$ and for any $\sigma\in G-\{1\}$, we have $\sigma(S_1)\neq S_1$. The Riemann-Roch Theorem ensures the existence of two positive integers $n<m$ such that the divisor $F=nS_0+mS_1$ is the polar divisor of a rational function $h\in K(X)$ (we can even set $n=1$). By applying Belyi's algorithm, we then obtain a polynomial $g\in K[T]$ of degree $d$ such that the morphism $f=g\circ h$ is a Belyi map. By construction, the divisor $dF$ is the polar divisor of $f$. We just have to check that the group $\aut_\qbar(f)$ is trivial. Any $\qbar$-automorphism $\sigma$ of $f$ induces an automorphism of $F$. In particular, it permutes the elements of the support of $S_0$, so that it belongs to $G$. Moreover, it also permutes the points of $S_1$ and thus $\sigma=1$, since $S_1$ was constructed in such a way that $S_1\neq\sigma(S_1)$ for any non-trivial element $\sigma\in G$.
\end{proof}

We can now define the {\bf relative Belyi degree} of $X$ as the minimal degree of an element of $\bel^0(X,K)$,
$$\deg_B(X,K)=\min\{\deg(f)\,\,|\,\,f\in\bel^0(X,K)\},$$
and Lemma~\ref{lemm1} directly leads to the relative version of Proposition~\ref{finite1}:

\begin{prop}\label{finite2} For any number field $K$, there exist finitely many $K$-isomorphism classes of curves of bounded relative Belyi degree.
\end{prop}

\begin{exem} As we have seen in Example 1.2, the $\qbar$-isomorphism class of elliptic curves having $j$-invariant $\frac{207646}{6561}$ has (absolute) Belyi degree $4$. An explicit model can be obtained by considering the elliptic curve $E$ defined by the affine Weierstrass equation
$$Y^2=(X+42)(X^2-42X+3033),$$
the cover being induced by the rational function
$$f=\frac1{8748}(X^2-6X+4545+6Y).$$
Since the $\qbar$-isomorphism group of $f$ is trivial, we have the identity
$$\deg_B(E,\bold Q)=4.$$
\end{exem}

\subsection{Primes of bad reduction, minimal regular models} Let $X$ be a curve of positive genus defined over a number field $K$. The inconvenient of the stable model of $X$ is the fact that it is generally defined over a non-trivial extension of $K$. In order to get rid of this problem, we replace it with the minimal regular model $\mathcal X$ of $X$, which is defined over the ring of integers $\mathcal O_K$ of $K$ (see for example~\cite{Chinburg}). We say that $p$ is a {\bf prime of bad reduction} if there exists a prime $\frak p$ of $\mathcal O_K$ lying above $p$ such that $\mathcal X$ has bad reduction at $\frak p$.

\subsection{A lower bound}

The following result is the analogue of Theorem~\ref{bound1}:

\begin{theo}\label{bound2} For any curve $X$ of positive genus defined over a number field $K$, we have the inequality
$$\deg_B(X,K)\geq p,$$
where $p$ is the greatest prime of bad reduction.
\end{theo}

\begin{proof} We proceed as in the proof of Theorem~\ref{bound1}, by proving that the minimal regular model of $X$ has good reduction at any prime $\frak p$ of $K$ not dividing the order of the monodromy group $G$ of $f$. The inequality of the theorem will again follow from the injection of $G$ in $S_n$, where $n$ is the degree of $f$. We can clearly work locally and replace $\mathcal O_K$ with its localization $\mathcal O_\frak p$ at $\frak p$. We therefore consider the minimal regular model of $X$ as being defined (by base change) over $\mathcal O_\frak p$. Following Proposition 5.3 in~\cite{Beckmann}, there exists a smooth model $\mathcal X_L$ of $X\otimes_KL$ defined over the localization $\mathcal O_\frak q$ of the ring of integers of a finite extension $L$ of $K$ at a prime $\frak q$ lying above $\frak p$. Following Proposition 3.1 in~\cite{Beckmann}, the model $\mathcal X_L$ descends to a smooth model $\mathcal X_N$ defined over an extension $N/K$ unramified at $\frak q_N=\frak q\cap N$. The main inconvenient is that the generic fiber of $\mathcal X_N$ will only be $L$-isomorphic (and generaly not $N$-isomorphic) to $X\otimes_KN$. Now, since $\aut_\qbar(f)$ is trivial, Lemma~\ref{lemm1} implies that $\mathcal X_N$ is actualy a model of $X\otimes_KN$. In particular the minimal regular model of $X\otimes_KN$ has good reduction at $\frak q_N$. We finally use Lemme 11 in~\cite{Liu}, which asserts that the minimal regular model of $X\otimes N$ is obtained  by (\'etale) base change from $\mathcal O_\frak p$ to $\mathcal O_{\frak q_N}$ of the minimal regular model of $X$. The result follows from the fact that smoothness descends under \'etale base change.\end{proof}

\begin{exem} For any prime number $p>3$, consider the elliptic curve $E_p$ over $\bold Q$ defined by the affine Weierstrass equation
$$Y^2=X^3+p$$
Since the $j$-invariant of this curve is equal to $0$, we have the identity $\deg_B(E_p,\qbar)=3$. Nevertheless, the above result gives the inequality $\deg_B(E_p,\bold Q)\geq p$. Indeed, $p^2$ is the greatest power of $p$ dividing the discriminant of $E_p$. This implies that the above equation is minimal at $p$ and, considered over $\bold Z_p$, it gives the N\'eron model of $E_p$ at $p$, which is therefore of type II (following Kodaira symbols, cf.~\cite{Silverman}, table 4.1) and has bad reduction at $p$. Since the N\'eron model is the smooth locus of the minimal regular model of $E_p$, the assertion follows.
\end{exem}

\bibliographystyle{alpha}
\bibliography{Belyidegree}

\end{document}